\documentclass[12pt,bezier]{conm-p-l}
\usepackage{amsfonts,latexsym,amscd}

\newtheorem{theorem}{Theorem}[section]

\newtheorem{proposition}[theorem]{Proposition}
\newtheorem{corollary}[theorem]{Corollary}
\newtheorem{lemma}[theorem]{Lemma}

\newcommand{\prop}[1]{Proposition \ref{#1}}
\newcommand{\lma}[1]{Lemma \ref{#1}}

\newcommand{\cor}[1]{Corollary \ref{#1}}

\newcommand{\rarrow}{\rightarrow}

\newcommand{\fd}{\operatorname{fd}}
\newcommand{\der}{\operatorname{Der}}

\newcommand{\chr}{\operatorname{char}}
\newcommand{\ho}{\operatorname{Ho}}

\newcommand{\spec}{\operatorname{Spec}}
\newcommand{\tor}{\operatorname{Tor}}
\newcommand{\Hom}{\operatorname{Hom}}
\newcommand{\ext}{\operatorname{Ext}}
\newcommand{\dtensor}{\operatorname{\otimes^{{\bf L}}}}

\newcommand{\calL}{{\mathcal L}}

\newcommand{\la}{\operatorname{\!\langle\!}}
\newcommand{\ra}{\operatorname{\!\rangle}}

\newcommand{\calM}{{\mathcal M}}

\newcommand{\calA}{{\mathcal A}}

\newcommand{\calO}{{\mathcal O}}
\newcommand{\calF}{{\mathcal F}}

\newcommand{\calC}{{\mathcal C}}

\begin{document}

\title[Characterizing Simplicial Algebras with Vanishing Homology]
{Characterizing Simplicial Commutative Algebras with Vanishing 
Andr\'e-Quillen Homology}
\author{James M. Turner}
\address{Department of Mathematics\\
Calvin College\\
3201 Burton Street, S.E.\\
Grand Rapids, MI 49546}
\email{jturner@calvin.edu}
\thanks{Partially supported by National Science Foundation (USA) grant 
DMS-0206647 and a Calvin Research Fellowship. He thanks the Lord for 
making his work possible.}
\date{\today}



\begin{abstract}
The use of homological and homotopical devices, such as $\tor$ and 
Andr\'e-Quillen homology, have found substantial use in 
characterizing commutative algebras. The primary category setting has 
been differentially graded algebras and modules, but recently 
simplicial categories have also proved to be useful settings. In this 
paper, we take this point of view up a notch by extending some  
recent uses of homological algebra in characterizing Noetherian 
commutative algebras to characterizing simplicial 
commutative algebras having {\it finite Noetherian homotopy} through 
the use of simplicial homotopy theory. These 
characterizations involve extending the notions of locally complete 
intersections and locally Gorenstein algebras to the simplicial 
homotopy setting.
\end{abstract}

\maketitle

\section*{Overview}

Following a program set forth by Grothendieck (see \cite{AF2}), a major research 
effort has been underway to characterize ring homomorphisms $f: R \to S$ 
of Noetherian rings. Use of homological devices have been pivotal to 
making such characterizations. For example, in \cite{Qui3} D. 
Quillen, in the process of developing a homology of commutative 
algebras, conjectured that the higner vanishing of this homology of 
$S$ over $R$, in the case $f$ is essentially of finite type and of 
finite flat dimension, charactertizes $f$ as a locally complete intersection 
homomorphism. Motivated by this conjecture, L. Avramov \cite{Avr} 
defined locally complete intersection homomorphisms more generally and 
established their properties, including providing a proof of Quillen's 
conjecture by a tour de force use of differentially graded and 
simplicial techniques. This fit into a larger program of Avramov with 
various collaborators to fulfill Grothendieck's program for rings and 
their homomorphisms. See, for example, \cite{Avr, AF, AFH}. In 
particular, characterizations of homomorphisms to be locally regular, 
complete intersection, Grothendieck, and Cohen-Macaulay were achieved.

Also motivated by Quillen's conjecture and a perspective of 
commutative algebra from a strictly simplicial viewpoint (see, for 
example, \cite{Goe, Mil}), a version of Quillen's conjecture was 
formulated and proved for simplicial commutative algebras with {\it 
finite Noetherian homotopy} \cite{Tur1,Tur2}, relying on techniques 
developed by analogy from the homotopy of spaces for the hootopy of 
simplicial commutative algebras. In the process, the notion of 
{\it homotopy complete intersection} was formulated and shown to be 
characterized by the vanishing of higher Andr\'e-Quillen homology. 
The drawback was that the extension was valid only when the $\pi_{0}$ 
has non-zero characteristic. Since such restrictions is not needed in 
the constant simplicial case, as the main result of \cite{Avr} 
clearly implies, then there is a gap in a full characterization of 
locally complete intersections through purely simplicial techniques.

The aim of this paper is twofold. The first is to begin to describe a 
program which extends the notions of complete intersection, 
Gorenstein, and Cohen-Macaulay to simplicial commutative algebras 
having Noetherian homotopy. This will involve summarizing the results 
of \cite{Tur1,Tur2} and extending some of the notions in 
\cite{AF,AFH}. In particular, we will define the notions of {\it homotopy 
Gorenstein} and {\it homotopy Cohen-Macaulay} algebras. 
The second aim of this paper will be to characterize, using these 
notions, simplicial algebras with finite Noetherian homotopy, finite 
flat dimensional homotopy, and finite Andr\'e-Quillen homology free of
conditions on the characteristic of $\pi_{0}$. Specifically, we will 
prove the following (see \S 2 for definitions):

\bigskip

\noindent {\bf Homology Characterization Theorem.} {\it Let $A$ be a 
simplicial commutative $R$-algebra ($R$ Noetherian) having finite Noetherian 
homotopy and $\fd_{R}(\pi_{*}A)<\infty$. If $D_{s}(A|R;-) = 0$ for $s\gg 0$, 
as a functor of $\pi_{0}A$-modules, then:
\begin{enumerate}
	\item[(a)] $A$ is locally homotopy Gorenstein;
	
	\item[(b)] If $\chr (\pi_{0}A) \neq 0$ then $A$ is a locally homotopy 
	complete intersection;
		
	\item[(c)] If $A$ is locally regularly 2-degenerate then $R \to 
	A$ is a locally complete intersection homomorphism; that is, for 
	each $\wp\in \spec (\pi_{0}A)$ there is a factorization (\ref{homfac}) such that 
	$\ker (\eta^{\prime})$ is generated by a regular sequence.
	
\end{enumerate}}

\bigskip

This paper is organized as folowed: \S 1 reviews the simplicial model 
structure for simplicial commutative algebras, as well as the tensor, 
hom, and differential structures and their derived functors. We also 
review relationships between certain simplicial categories and 
differentially graded categories. 
\S 2 reviews the definition of homotopy complete intersections and 
introduces the notions of homotopy Gorenstein and homotopy 
Cohen-Macaulay algebras. We end with a proof of the Homology 
Characterization Theorem.

\bigskip

\noindent {\bf Acknowledgements:} The author would like to thank 
Lucho Avramov for several discussions over the years which 
greatly assisted him in the development of this paper's content and 
for giving him access to \cite{AFHal}. He would also like to thank
Brooke Shipley for alerting him to the paper \cite{SS}.

\bigskip 

\section{Homotopy theory of simplicial commutative algebras and 
simplicial modules}

\subsection{Simplicial model category structures}

Fix a commutative ring $R$ with unit and let $\calA$ denote the 
category of commutative $R$-algebras with unit. If $B\in \calA$, let 
$\calA_{B}$ denote the subcategory of objects $(A,\epsilon)$ with 
$\epsilon : A \to B$ an $R$-algebra map. Given any 
category $\calC$, we let $s\calC$ denote the category of simplicial 
objects over $\calC$. Finally, for $A\in s\calA$, let $\calM_{A}$ 
denote the category of simplicial $A$-modules. 

Let $f: A \to B$ be a map in $s\calA$. Recall from \cite{Qui1} that 
$f$ is defined to be a:

\begin{enumerate}
	
	\item {\it weak equivalence} provided $f$ is a weak equivalence of 
	simplicial groups;
	
	\item {\it fibration} provided $f$ is a fibration of simplicial 
	groups;
	
	\item {\it cofibration} if and only if $f$ is a retract of an {\it almost 
	free map}, i.e. a map $f^{\prime}: A^{\prime} \to B^{\prime}$ such that 
	$f^{\prime}$ makes $B^{\prime}$ free as an almost simplicial 
	$A^{\prime}$-algebra (that is, without $d_{0}$).
	
\end{enumerate}	

By \cite[\S II.3]{Qui1}, this structure makes $s\calA_{S}$ a closed 
model category for any fixed $S\in s\calA$. 

Fixing $A\in s\calA$, let $f: K \to L$ be a map in $\calM_{A}$.
Following \cite{Qui1, SS}, we will say that $f$ is defined to be a:

\begin{enumerate}
	
	\item {\it weak equivalence} provided $Nf$ is a weak equivalence of 
	chain complexes;
	
	\item {\it fibration} provided $Nf$ is a level-wise surjection in 
	positive degrees;
	
	\item {\it cofibration} if and only if $Nf$ is a level-wise 
	monomorphism whose cokernel in each degree $k\geq 0$ is a projective 
	$N_{k}A$-module. 
	
\end{enumerate}	

Here $N: \calM_{A}\to Mod_{NA}$ is the normalized chain functor. 
By \cite[\S II.3]{Qui1}, this structure makes $\calM_{A}$ into a closed 
model category. Furthermore, part of the main theorem of \cite{SS} states

\bigskip

\noindent {\bf Schwede-Shipley Theorem:} {\it The Dold-Kan 
correspondence $$N: \calM_{A} \Longleftrightarrow Mod^{+}_{NA}: K$$ is a 
Quillen equivalence of symmetric monoidal closed model categories.}

\bigskip

\subsection{Tensor and Tor-modules}

For $A\in s\calA$, level-wise tensor product gives a functor 
$$\otimes_{A}:\calM_{A}\times\calM_{A} \to \calM_{A}.$$
In turn, this tensor product
induces the {\it derived tensor product} $$\otimes_{R}^{{\bf L}}: \ho 
(\calM_{A}) \times \ho (\calM_{A}) \to \ho (\calM_{A})$$ on the 
homotopy category, which is
defined by $$K\otimes_{A}^{{\bf L}}L := X\otimes_{A}Y,$$ where $X$ and 
$Y$ are cofibrant replacements for $K$ and $L$ in $\calM_{A}$, respectively. 

Note also that $\otimes_{A}$ induces a functor on $s\calA$ which 
descends to a functor $\otimes_{A}^{{\bf L}}$ on $s\calA_{B}$, for any 
fixed $B$.

Next, define the {\it Tor-modules} for $K,L \in \calM_{A}$ to be
$$
\tor_{s}^{A}(K,L) := \pi_{s}(K\otimes_{A}^{{\bf L}}L), \quad s\geq 0.
$$
For each such $s$, $\tor_{s}^{A}$ is a $\pi_{0}A$-module.

A key method for computing $\tor_{*}^{A}$ is the following device found in 
\cite[II.6]{Qui1}. 

\bigskip

\noindent {\bf Kunneth Spectral Sequence:} {\it For $K,L\in \calM_{A}$,
there is a first quadrant spectral sequence 
$$E^{2}_{*,*}=\tor_{*}^{\pi_{*}A}(\pi_{*}K,\pi_{*}L) \Longrightarrow 
\tor_{*}^{A}(K,L).$$}

\bigskip

\subsection{Ext-modules}

As noted above, in \cite{SS} it is shown that the normalization 
functor induces a functor $N : \calM_{A} \to Mod_{NA}^{+}$ which is an 
equivalence of categories. We will therefore use $N$ to define the 
{\it Ext-modules} for $K,L \in \calM_{A}$. First, if $T$ is a 
(commutative) DG ring and $U,V$ are DG T-modules, let
$$
\ext_{T}^{s}(U,V) := H_{-s}(\Hom_{T}(P_{*},I_{*}))
$$
where $P_{*}\to U$ is a DG projective replacement and $V \to I_{*}$ is a DG injective
replacement for $NL$ in $NA-Mod$. Here the {\it Hom-complex} is defined 
as in \cite[\S 1]{AFL}. We then define, for $K,L\in \calM_{A}$, $$\ext_{A}^{s}(K,L) := 
\ext_{NA}^{s}(NK,NL).$$ 
A basic property of $\ext$ is then given by:

\begin{lemma}\label{exteq}
	\cite[(1.8)]{AFL}
	If $A \stackrel{\sim}{\rarrow} B$ is a weak equivalence of 
	simplicial algebras augmented over a field $\ell$, then 
	$$\ext_{A}^{*}(\ell, A)\cong \ext_{B}^{*}(\ell,B).$$
\end{lemma}	

\bigskip 

Now assume that $\ell$ is a field of characteristic 0. Let $\Re_{\ell}$ denote the 
category of commutative $\ell$-algebras over $\ell$, henceforth 
refered to as {\it rational $\ell$-algebras}. Let $s_{+}\Re_{\ell}$ and 
$ch_{+}\Re_{\ell}$ denote the category of connected simplicial rational 
$\ell$-algebras and the category of connected differentially graded 
rational $\ell$-algebras, respectively.

\bigskip

\noindent {\bf Quillen's Theorem.} \cite[p. 223]{Qui2} \; {\it Normalization 
induces a functor $N : s_{+}\Re_{\ell} \to ch_{+}\Re_{\ell}$ which is an equivalence of 
closed model categories.}

\bigskip

Our aim now is to prove the following variation of \cite[(4.1)]{AFL}.

\begin{proposition}\label{extprop}
	Let $\psi: \ell[x] \to T$ be a map of augmented DG rings over 
	$\ell$ from a free commutative graded $\ell$-algebra on one generator $x$, $|x| = n >0$,
	and let $M$ be a DG $\ell[x]$-module. Assume furthermore 
	that
	\begin{enumerate}
		\item $H_{0}(T)$ is Noetherian and each $H_{i}(T)$ is a finitely 
		generated $H_{0}(T)$-module for $i\in \mathbb Z$;
		\item $H(M)$ is bounded above.
	\end{enumerate}
	Letting $F(\psi) := T\otimes_{\ell[x]}^{{\bf L}}\ell$, there is 
	an isomorphism of graded $\ell$-modules $$\ext_{T}^{*}(\ell, 
	M\dtensor_{\ell[x]}T) \cong 
	\ext_{\ell[x]}^{*}(\ell,M)\otimes_{\ell} 
	\ext^{*}_{F(\psi)}(\ell,F(\psi)).$$
\end{proposition}

\noindent {\it Proof.} By inspection of the statement of 
\cite[4.1]{AFL}, the key difference involves replacing an augmented 
DG ring $R$, for which $H_{*}(R) \cong H_{0}(R)$ is Noetherian, with 
$N(S_{\ell}(n))$. In their proof of (4.1), the condition on $R$ 
insures that:
\begin{enumerate}
	\item There is a factorization $$R\rarrow X \stackrel{\simeq}{\rarrow} \ell$$ with $X$ a free $R$-algebra such 
     that $X\stackrel{\simeq}{\rarrow} X^{\prime}$ with $X^{\prime}$ 
     finite type over $R$ with bounded below generators. 
	
	 \item There is a DG injective $R$-resolution $M \to I$ with $I$ 
	 bounded above.
\end{enumerate}
(1) and (2) insures that \cite[(1.10)]{AFL} can be applied to establish \cite[(4.9)]{AFL}. 

To show (1) in our context, let $\ell[x,y]$ be the free DG $\ell$-algebra 
such that $|y| = n+1$ and $\partial y = x$. Then there is a factorization 
$$
\ell[x] \rarrow \ell[x,y] \stackrel{\sim}{\rarrow} \ell
$$
with the required properties. Finally, (2) can be found in 
\cite[(9.3.2.1)]{AFHal}.
\hfill $\Box$  

\bigskip

Let $S_{\ell}(n)$ be the free commutative $\ell$-algebra 
generated by the Eilenberg-MacLane object $K(\ell,n)$. 
Let $A\in s\Re_{\ell}$ and let $\phi : S_{\ell}(n) \to A$ be a map 
in $s\Re_{\ell}$. $\phi$ determines a cofibration sequence in $\ho(s\Re_{\ell})$
$$
S_{\ell}(n) \stackrel{\phi}{\rarrow} A \rarrow \calF_{\phi}
$$
where $\calF_{\phi} := A\otimes_{S_{\ell}(n)}^{{\bf L}}\ell$. 

\bigskip 

\begin{corollary}\label{extfac}
	If $A\in s_{+}\Re_{\ell}$ with $\pi_{*}A$ 
	a finite graded $\ell$-module, then $$\ext_{A}^{*}(\ell, A) \cong 
	\ext_{S_{\ell}(n)}^{*}(\ell,S_{\ell}(n))\otimes_{\ell} 
	\ext^{*}_{\calF_{\phi}}(\ell,\calF_{\phi}).$$ 
\end{corollary}

\noindent {\it Proof.} 
Note first that $N(S_{\ell}(n)) \simeq \ell[x]$ with $|x|=n$. Since 
the result follows immediately from \prop{extprop} for n odd 
($\ell[x]$ is exterior on $x$), we assume that n is even. 
Using Quillen's Theorem, let $L(k)\in s\Re_{\ell}$ satisfy $NL(k) \simeq
\ell[x]/(x^{k})$. Thus, by \prop{extprop}:
$$
\ext_{A}^{*}(\ell,L(k)\dtensor_{S_{\ell}(n)}A) \cong 
\ext_{S_{\ell}(n)}^{*}(\ell,L(k))\otimes_{\ell} \ext^{*}_{\calF_{\phi}}(\ell,\calF_{\phi}).
$$
(Note that the Schwede-Shipley Theorem implies that $N\calF_{\phi}\simeq 
F(N(\phi))$.)

Now $NS_{\ell}(n) \to NL(k)$ is equivalent to a 
map that is an isomomorphism in degrees $<nk$. 
By a Kunneth spectral sequence argument, since $\pi_{*}A$ is bounded the map $A\to 
L(k)\dtensor_{S_{\ell}(n)}A$ induces a $\pi_{*}$-injection which 
is a $\pi_{*}$-isomorphism through degree $nk$, for $k\gg 0$. In 
particular, as $\pi_{*}A$-modules: 
$$\tor_{*}^{S_{\ell}(n)}(L(k),A) \cong \pi_{*}A\times 
\Sigma^{nk+1}\pi_{*}A.$$ Thus 
the map $A\to \lim_{k}(L(k)\dtensor_{S_{\ell}(n)}A)$ is a weak 
equivalence, whose normalization is equivalent to an isomorphism. 
Furthermore, using Quillen's Theorem and the Schwede-Shipley Theorem, the induced map 
$NA \to (\ell[x]/(x^{k}))\dtensor_{\ell[x]}NA$ 
can be shown to be equivalent to a split injection, for $k\gg 
0$, as $NA$-modules. Thus the result now follows from an argument using 
\lma{exteq} and the Milnor sequence \cite[\S 3.5]{Wei}. \hfill $\Box$

\bigskip 

\noindent {\bf Note:} \cor{extfac} holds over fields $\ell$ of 
arbitrary characteristic when $n = 1$.

\bigskip

\subsection{Differentials and Andr\'e-Quillen homology}

Let $A$ be a commutative ring and $B$ a commutative $A$-algebra. If 
$M$ is a $B$-module, recall that an $A$-module map $f: B \to M$ is a 
{\it derivation} provided $f(xy) = xf(y) + yf(x)$. Let 
$\der_{A}(B,M)$ be the $A$-module of derivations. The functor $M \mapsto 
\der_{A}(B,M)$ is representable: there is a canonically 
defined $B$-module $\Omega_{B|A}$, called the {\it differentials} of 
$B$ over $A$, such that there is a natural isomorphism: 
$$
\der_{A}(B,M) \cong \Hom_{B}(\Omega_{B|A},M).
$$

From the differentials, the {\it cotangent complex} of a simplicial 
commutative $R$-algebra $A$ is defined by
$$
\calL(A|R) := \Omega_{X|R}\otimes_{X}A,
$$
where $X \stackrel{\sim}{\rarrow} A$ is a cofibrant replacement of 
$A$ in $\ho(s\calA)$. The {\it Andr\'e-Quillen homology} of $A$ over $R$ 
with coefficients in the $A$-module $M$ is then defined to be
$$
D_{*}(A|R;M) := \pi_{*}(\calL(A|R)\otimes_{A}M).
$$
We now recall two important properties of Andr\'e-Quillen homology:
\begin{enumerate}
	\item (Transitivity Sequence) Given maps $A \to B \to C$ in $s\calA$ 
	and a $C$-module $M$, there is a long exact sequence:
	$$
	\ldots \to D_{s+1}(C|B;M) \to D_{s}(B|A;M) \to D_{s}(C|A;M) \to \to 
	D_{s}(C|B;M) \to \ldots;
	$$
	\item (Flat Base Change) For $A,B \in s\calA$ and $M$ an 
	$A\dtensor_{R}B$-module
	$$
	D_{*}(A\dtensor_{R}B|B;M) \cong D_{*}(A|R;M).
	$$
\end{enumerate}
A further useful relationship between homotopy and Andr\'e-Quillen 
homology is given by the following:

\bigskip

\noindent {\bf Hurewicz Theorem:} {\it For a connected simplicial 
$R$-algebra $A$ over a field
$\ell$, there is a homomorphism $$h_{*}: \tor_{*}^{R}(A,\ell) \to 
D_{*}(A|R;\ell)$$ for which $h_{*}$ is an isomorphism in degrees $\leq 
n$ provided $A$ is (n-1)-connected.}

\bigskip

Let $\epsilon: A \to \ell$ be a commutative algebra over a field. 
Let $I = \ker (\epsilon)$. Define the {\it indecomposables} of $A$ to 
be the $\ell$-module $QA :=I/I^{2}$. A well known result \cite{Goe,Mil} 
for $A$ a supplemented $\ell$-algebra is
$$
\Omega_{A|\ell}\otimes_{A}\ell \cong QA.
$$
We thus define the Andr\'e-Quillen homology of a simplicial 
supplemented $\ell$-algebra $A$ by
$$
H^{Q}_{*}(A) := \pi_{*}QX,
$$
where $X \stackrel{\sim}{\rarrow} A$ is a cofibrant resolution of $A$ 
as simplicial supplemented algebras. Thus we have
$$
H^{Q}_{*}(A) \cong D_{*}(A|\ell;\ell).
$$

We now return to inspecting maps $\phi : S_{\ell}(n) \to 
A$ of simplicial supplemented $\ell$-algebra. 
Our aim is to prove

\begin{proposition}\label{aqprop}
	Let $\phi: S_{\ell}(n) \to A$ be a map of (n-1)-connected simplicial 
	supplemented $\ell$-algebras such that $H^{Q}_{n}(\phi)$ is an 
	injection. If $\pi_{*}A$ is a finite graded $\ell$-module and 
	$\pi_{*}(\calF_{\phi})$ is unbounded then $H^{Q}_{*}A$ is unbounded.
\end{proposition}

\noindent (Recall that a positively graded module $M$ is {\it unbounded} provided 
$M_{t}\neq 0$ for infinitely many $t$.)

\bigskip 

To prove this proposition, we adapt the proof of \cite[(4.2)]{AH}. Assume $\chr \ell = 0$.
Recall for a commutative DG $\ell$-algebra $T$ that a {\it minimal 
model} for $T$ is a free DG-algebra $(\ell[X],\partial)$ such that 
$\partial X \subseteq I^{2}$, where $I$ is the augmentation ideal of $\ell[X]$, 
together with a weak equivalence $\ell[X] \stackrel{\sim}{\rarrow} T$. 
For existence and properties of minimal models, see \cite{BG}.
Note that from Quillen's Theorem (see also \cite[Thm. 9.5]{Qui3}, if $A$ is a simplicial 
supplemented $\ell$-algebra and $\ell[X]$ is a minimal model for $NA$,
then 
$$
H^{Q}_{*}(A) \cong QX.
$$

\bigskip

\noindent {\it Proof of \prop{aqprop}.} First, assume that $\chr \ell 
\neq 0$. If $\pi_{*}A$ and $H^{Q}_{*}A$ are both finite graded 
$\ell$-modules, then, by the Algebraic Serre Theorem \cite{Tur1}, 
$A\cong\bigotimes_{I}S_{\ell}(1)$ in the homotopy category. 
Thus $n = 1$ and $\calF_{\phi}\cong\bigotimes_{I}S_{\ell}(1)$ with 
$|J| = |I|-1$. It follows that $\pi_{*}\calF_{\phi}$ is bounded.

Now assume that $\chr \ell = 0$. Let $(\ell[X],\partial)$ be a minimal 
model for $NA$ and let $\ell[x]$ be a minimal model for 
$N(S_{\ell}(n))$. By the assumption that $H^{Q}_{n}(\phi)$ is 
injective, it follows that we may assume that $x\in X$ and that the 
map $\ell[x]\to \ell[X]$ induced by the inclusion $\{x\}\subseteq X$ 
is equivalent to $N(\phi)$. It follows from Quillen's Theorem that if 
$X = \{x\}\bigcup Y$, then $\ell[Y]$ is a minimal model for $N\calF_{\phi}$. 
Thus writing $\ell[X] \cong \ell[x]\otimes\ell[Y]$, we can express, for 
$u\in \ell[Y]$,
$$
\partial (1\otimes u) = 1\otimes\bar\partial u + x\otimes\calO u.
$$
It follows that $\calO : (\ell[Y],\bar\partial) \to 
(\ell[Y],\bar\partial)$ is a derivation of 
degree -n-1. 

Let $J$ be the augmentation ideal of $\ell[Y]$. Let $\calO = 
\Sigma_{i\geq 1}\calO_{i}$ with $\calO_{i}(J)\subseteq J^{i}$. Thus 
$\calO_{1}: J \to J$ is a derivation and, hence, induces
$$
\calO_{1} : J/J^{2} \rarrow J/J^{2}.
$$

\bigskip

\noindent {\bf Claim:} {\it There exists an element $u$ in the $\ell$-dual 
$(J/J^{2})^{*}$ such that $(\calO_{1}^{*})^{n}u\neq 0$ for all $n\geq 0$.}

\bigskip

\noindent It follows from this claim that $H^{Q}_{*}(\calF_{\phi})$ and, 
hence, $H^{Q}_{*}(A)$ are unbounded.   

\bigskip 

To establish the claim, assume for each $u\in (J/J^{2})^{*}$ there is 
an $n\gg 1$ such that $(\calO_{1}^{*})^{n}u = 0$. Following the 
same argument as in \cite[p.181]{AH}, this implies that, for each $u\in 
J^{*}$, $(\calO^{*})^{n}u = 0$ for $n\gg 1$.

Now, consider the exact sequence $$0\to J \stackrel{\tau}{\rarrow} I \to J \to 0,$$ where 
$I$ is the augmentation ideal of $\ell[X]$ and $\tau (w) = z\otimes w$.  
Dualizing and applying cohomology, we obtain a long exact sequence
$$
\ldots \to H^{i-1}(J^{*}) \stackrel{\delta}{\rarrow} H^{i+n}(J^{*}) 
\to (H_{i+n}(I))^{*} \stackrel{\tau^{*}}{\rarrow} H^{i}(J^{*}) \to \ldots
$$
An easy computation shows that $\delta = H(\calO^{*})$. By our 
assumption on the finiteness of $\pi_{*}A$, $H_{i}(I) = 0$ for $i\geq N$, 
$N\gg 0$. Thus $\delta$ is injective for $i\geq N$. Since we 
are assuming $\pi_{*}\calF_{\phi}$ is unbounded, this contradicts our 
local nilpotency condition on $(\calO)^{*}$. Thus our claim is 
established. \hfill $\Box$

\section{Characterizing Simplicial Commutative Algebras}

We focus on extending the characterizations of homomorphisms $R \to 
S$ of Noetherian rings achieved in \cite{Avr,AF,AFH} to simplicial 
commutative $R$-algebras. To set in what direction this extension is 
to take, we view $S$ as a constant simplicial $R$-algebra. To achieve 
a suitable type of extension, the notion of Noetherian needs to be 
spelled out for simplicial algebras. Such a notion was delineated and 
explored in \cite{Tur2}, motivated by an analogous for concept for DG 
rings described in \cite{AF}, which we now describe.

A simplicial commutative algebra $A$ is said to have {\it Noetherian 
homotopy} provided
\begin{enumerate}
	
	\item $\pi_{0}A$ is a Noetherian ring;
	
	\item each $\pi_{m}A$ is a finite $\pi_{0}A$-module.

\end{enumerate}
We furthermore say that $A$ has {\it finite Noetherian homotopy} 
provided that (2) is replaced by
\begin{enumerate}
	
	\item[$(2)^{\prime}$] $\pi_{*}A$ is a finite graded $\pi_{0}A$-module. 

\end{enumerate}

The key to characterizing simplicial commutative $R$-algebras with 
Noetherian homotopy through homotopical/homological methods is to 
locally reduce to connected simplicial algebras over a field. Such 
objects yield more information under homological scrutinity. This 
approach was pioneered by L. Avramov \cite{Avr1} through the notion 
of DG fibre and used with great effect in \cite{AF,Avr}. To adopt this 
approach in the simplicial setting, the following extension of the 
main theorem of \cite{AFH} is needed

\bigskip

\noindent {\bf Factorization Theorem.} \cite[(2.8)]{Tur2} {\it 
	Suppose $A$ is a simplicial commutative $R$-algebra with Noetherian 
	homotopy, $R$ a Noetherian ring, and $\wp\in \spec(\pi_{0}A)$. 
	Then there is a simplicial commutative algebra $A^{\prime}$ with 
	Noetherian homotopy, such that $\pi_{*}A^{\prime} \cong \widehat{\pi_{*}A}$, and
	there exists a (complete local) Noetherian $R^{\prime}$ that fits into the 
	following commutative diagram in $\ho (s\calA_{k(\wp)})$
\begin{equation}\label{homfac}
\begin{array}{ccc}
R & \stackrel{\eta}{\longrightarrow} &
A \\[1mm]
\phi \downarrow \hspace*{10pt}
&&
\hspace*{10pt} \downarrow \psi \\[1mm]
R^{\prime} & \stackrel{\eta^{\prime}}{\longrightarrow} & A^{\prime}
\end{array}
\end{equation}
     with the following properties:
	 \begin{enumerate}
		 
		 \item $\phi$ is a flat map and its closed fibre $R^{\prime}/\wp 
		 R^{\prime}$ is regular;
		 
		 \item $\psi$ is a flat $D_{*}(-|R;k(\wp))$-isomomorphism;
		 
		 \item $\eta^{\prime}$ induces a surjection $\eta^{\prime}_{*}: R^{\prime} \to 
		 \pi_{0}A^{\prime}$;
		 
		 \item $\fd_{R} (\pi_{*} A)$ finite implies that $\fd_{R^{\prime}} 
		 (\pi_{*} A^{\prime})$ is finite
		 
	\end{enumerate}	 }

\bigskip 

\noindent We will call a choice of diagram (\ref{homfac}) a {\it 
factorization} for a simplicial $R$-algebra $A$ with Noetherian 
homotopy. Also, recall that, for an $R$-module $M$, $\fd_{R}M$ is 
the {\it flat dimension} of $M$.

\bigskip 

We now describe extensions of the notions of locally complete 
intersection, locally Gorenstein, and locally CM (i.e. Cohen-Macaulay) for homomorphisms 
of Noetherian rings, as described in \cite{Avr,AF,AFH}, to simplicial 
algebras with Noetherian homotopy.

To begin, let $A$ be a connected simplicial supplemented $\ell$-algebra, $\ell$ 
a field. We then declare that $A$ is
\begin{enumerate}
	
	\item {\it a homotopy complete intersection} provided there is a 
	finite set $I$ with $A\cong\bigotimes_{I}S_{\ell}(1)$ in $\ho (s\calA)$;
	
	\item {\it homotopy CM} provided there 
	exists an $n\in \mathbb Z$ such that $\ext^{i}_{A}(\ell,A) = 0$ for $i\neq n$;
	
	\item {\it homotopy Gorenstein} provided $A$ is homotopy CM and 
	$\dim_{\ell}\ext_{A}^{n}(\ell,A) = 1$.
	
\end{enumerate}

To get a sense of how these notions fit together and apply to basic examples of 
simplicial supplemented algebras, we prove the following:

\begin{proposition}\label{basic}
	\begin{enumerate}
		
		\item[(a)] $S_{\ell}(n)$ is a homotopy complete intersection when $n =1$ 
		and homotopy Gorenstein in general when $\chr\ell = 0$.
		
        \item[(b)] We have the string of implications $$\text{homotopy complete 
		intersection}\Longrightarrow \text{homotopy Gorenstein}\Longrightarrow 
		\text{homotopy CM}.$$
				
	\end{enumerate}
\end{proposition}	

\noindent {\it Proof.} To prove (a), we note that $NS_{\ell}(n) 
\simeq \ell[x]$, a free commutative graded algebra, on one generator 
$x$ with $|x|=n$, and zero differential. Since $\ell[x]$ is Gorenstein it is therefore 
homotopy Gorenstein (see \cite[(18.1) \& (21.3)]{Mat}). The result now follows from 
\cite[(1.8)]{AFL}.

To prove (b), it is clearly enough to prove 
that a homotopy complete intersection is homotopy Gorenstein. But for 
this, note that an inclusion $\phi: S_{\ell}(1) \hookrightarrow 
\bigotimes_{I}S_{\ell}(1)$ onto one factor has homotopy fibre 
$\calF_{\phi}\cong \bigotimes_{J}S_{\ell}(1)$ with $|J| = |I|-1$. The 
result now follows from an induction using (a) and \prop{extfac}.
\hfill $\Box$

\bigskip
Consider now a simplicial commutative $R$-algebra $A$ over a field 
$\ell$ such that the unit map $R\to \pi_{0}A$ is surjective. We will then 
say that $A$ is a {\it homotopy complete intersection} (resp. {\it 
homotopy Gorenstein}, {\it homotopy CM}) {\it over $\ell$} provided 
$A\dtensor_{R}\ell$ is a homotopy complete intersection (resp. 
homotopy Gorenstein, homotopy CM). 

Finally, for a Noetherian ring $R$ and a simplicial commutative $R$-algebra 
$A$ with Noetherian homotopy, we say that $A$ is a {\it locally 
homotopy complete intersection} (resp. {\it locally homotopy 
Gorenstein}, {\it locally homotopy CM}) provided that for each $\wp\in 
\spec (\pi_{0}A)$ there is a factorization (\ref{homfac}) such that the simplicial 
$R^{\prime}$-algebra $A^{\prime}$ is a homotopy complete intersection 
(resp. homotopy Gorenstein, homotopy CM) over the residue field 
$k(\wp)$. 

\bigskip 

Now we proceed to proving our main result. Before we do so, we need a 
technical definition and result. First, given a simplicial $R$-algebra $A$ 
over a field $\ell$, with $(R,m)$ local Noetherian and $\eta: R \to \pi_{0}A$ a 
surjection, let $x_{1},\ldots,x_{n}\in m$ be a maximal regular 
sequence in $\ker (\eta)$ which extends to a minimal generating set 
for $\ker (\eta)$. We then say that $A$ is {\it regularly r-degenerate} 
provided the Kunneth spectral sequence
$$
\tor_{*}^{R/(x_{1},\ldots,x_{n})}(\pi_{*}A,\ell) \Longrightarrow 
\tor_{*}^{R/(x_{1},\ldots,x_{n})}(A,\ell)
$$
degenerates at the $E^{r}$-term. For a general simplicial commutative 
$R$-algebra $A$ with Noetherian homotopy, we declare that $A$ is {\it 
locally regularly r-degenerate} if for each $\wp \in \spec 
(\pi_{0}A)$, there is a factorization such that the simplicial 
$R^{\prime}$-algebra $A^{\prime}$ is regularly r-degenerate over 
$k(\wp)$.  

We will also need the following result.

\begin{lemma}\label{seq}
	For a simplicial $R$-algebra $A$ over a field $\ell$, let $x\in 
	\ker (\eta)$ be a non-zero divisor in $R$. Then the sequence 
	$$R/(x)\dtensor_{R}\ell \to A\dtensor_{R}\ell \to A\dtensor_{R/(x)}\ell$$ 
	is a cofibration sequence of simplicial supplemented $\ell$-algebras.
\end{lemma}

\noindent {\it Proof.} It is enough to check that 
$$(A\dtensor_{R}\ell)\dtensor_{(R\dtensor_{R/(x)}\ell)}\ell \simeq 
A\dtensor_{R/(x)}\ell$$ but this is a straightforward computation.
\hfill $\Box$

\bigskip

We now have reached the main goal of this paper. 

\bigskip 

\noindent {\it Proof of the Homology Characterization Theorem.} 
(a) For $\wp\in\spec (\pi_{0}A)$, choose a 
factorization (\ref{homfac}) of $A$. Under the hypotheses on $A$, if 
$\chr k(\wp)\neq 0$ then $A^{\prime}\dtensor_{R^{\prime}}\ell$ is a homotopy 
complete intersection, by Theorem B of \cite{Tur2}. Thus it is homotopy 
Gorenstein by \prop{basic} (b). So assume that $\chr k(\wp) = 0$.

Since $X := A^{\prime}\dtensor_{R^{\prime}}\ell$ is connected, we may 
assume that $X$ is (n-1)-connected for some $n\geq 1$. Let $\phi: 
S_{\ell}(n) \to X$ be a map so that $H^{Q}_{n}(\phi)(x)$ is a 
basis member of $H^{Q}_{n}(X)\cong D_{n}(A|R;k(\wp))$ (by Flat Base 
Change and the Factorization Theorem). Here we write 
$H^{Q}_{n}(S_{\ell}(n))\cong \ell\la x\ra$.

Now, by Flat Base Change and the Transitivity Sequence, 
$$
H^{Q}_{s}(\calF_{\phi}) \cong
\begin{cases}
	H^{Q}_{s}(X) & \quad s \neq n; \\
    H^{Q}_{n}(X)/\ell\la x \ra & \quad s = n.
\end{cases}
$$
Also, from the finite flat dimension condition and the Kunneth 
spectral sequence, $\pi_{*}X$ is a finite graded $\ell$-module. 
It follows from \prop{aqprop} that $\pi_{*}\calF_{\phi}$ is a finite graded 
$\ell$-module. Thus, by an induction, using \prop{basic} (a), we may 
assume that $\calF_{\phi}$ is homotopy Gorenstein. But now combining 
\prop{basic} (a) with \cor{extfac} it follows that $X$ is homotopy 
Gorenstein.

\bigskip 

\noindent (b) This is just Theorem B of \cite{Tur2}.

\bigskip 

\noindent (c) We adapt another argument of L. Avramov and S. Halperin \cite[(4.1)]{AH}.  
Again, choose a factorization (\ref{homfac}) for $A$ at $\wp\in \spec 
(\pi_{0}A)$. For the unit map $\eta^{\prime}: (R^{\prime},m^{\prime})\to 
A^{\prime}$, let $x_{1},\ldots,x_{n}\in m^{\prime}$ be a 
maximally regular sequence which extends to a minimal generating set 
for $\ker (\eta^{\prime})$. Note that if $x_{1},\ldots,x_{n}$ fails 
to generate this kernel, then, by the main result of \cite{AB}, 
$\tor^{R^{\prime}/(x_{1},\ldots,x_{n})}_{*}(\pi_{*}A^{\prime},k(\wp))$ 
is unbounded. Since we are assuming $A$ is 
locally regularly 2-degenerate, we may assume that 
$\tor^{R^{\prime}/(x_{1},\ldots,x_{n})}_{*}(A^{\prime},k(\wp))$ 
is unbounded. 

Now, for $x\in \ker (\eta^{\prime})$ a non-zero divisor in 
$R^{\prime}$, we have 
$$
H^{Q}_{*}(R^{\prime}/(x)\dtensor_{R^{\prime}}\ell) \cong 
D_{*}(R^{\prime}/(x)|R^{\prime};k(\wp)) \cong k(\wp)\la u |\; |u| =1\ra,
$$
by \cite[(6.25)]{And}. Thus $R^{\prime}/(x)\dtensor_{R^{\prime}}k(\wp) 
\cong S_{k(\wp)}(1)$ in the homotopy category, by \cite[(2.1.3)]{Tur1}. Thus, by an induction, 
using \lma{seq}, there is a cofibration sequence
$$
S_{\ell}(1) \rarrow 
A^{\prime}\dtensor_{(R^{\prime}/(x_{1},\ldots,x_{s-1}))}k(\wp) \rarrow 
A^{\prime}\dtensor_{(R^{\prime}/(x_{1},\ldots,x_{s}))}k(\wp)
$$
for each $1\leq s \leq n$. By the Factorization Theorem and our discussion above, 
there exists an $1\leq s \leq n$ such that
\begin{enumerate}
	\item $\tor_{*}^{R^{\prime}/(x_{1},\ldots,x_{s-1})}(A,k(\wp))$ 
	is bounded;
	
	\item $\tor_{*}^{R^{\prime}/(x_{1},\ldots,x_{s})}(A,k(\wp))$ 
	is unbounded.
\end{enumerate} 	
Furthermore, by Flat Base Change, the Transitivity Sequence, and our 
finiteness assumptions, we also have
\begin{enumerate}
	\item[(3)] $D_{*}(A^{\prime}|R^{\prime}/(x_{1},\ldots,x_{t});k(\wp))$ 
     is bounded for all $1\leq t \leq n$.  
\end{enumerate}
Applying \prop{aqprop}, (1) and (2) together imply that 
$D_{*}(A^{\prime}|R^{\prime}/(x_{1},\ldots,x_{s-1});k(\wp))$ 
is unbounded, contradicting (3). We therefore conclude that 
$\ker(\eta^{\prime}) = (x_{1},\ldots,x_{n})$. 
\hfill $\Box$

\bigskip

\noindent {\bf Remark:} \prop{basic} and the Homology Characterization Theorem shows 
that the standard stratified characterization of Noetherian rings 
$$
\text{Complete Intersection} \subseteq \text{Gorenstein}
$$
extends to our present homotopy setting, yet sensitivity to 
differences in characteristic appear. Nevertheless, we now
begin to fill a gap between the main result of \cite{Avr}, which is independent of 
characteristic, and the main results of \cite{Tur1,Tur2}, which are 
valid only in non-zero characteristics. In particular, we have the 
following consequence of part (c) of this theorem.

\bigskip

\noindent {\bf Corollary} {\it Let $\phi: R \to S$ be a homomorphism of 
Noetherian rings with $\fd_{R}S< \infty$. Then $D_{s}(S|R;-) = 
0$ for $s\gg 0$ implies that $\phi$ is a locally complete intersection
homomorphism.}


\bigskip 

\begin{thebibliography}{99}

\bibitem{And}
M. Andr\'e, {\em Homologie des alg\`ebres commutatives}, Die Grundlehren
der Mathematischen Wissenschaften 206, Springer-Verlag, 1974.

\bibitem{AB}
M. Auslander and D. Buchsbaum, ``Codimension and multiplicity,'' 
Annals of Math. (3) 68 (1958), 625-657. 

\bibitem{Avr1}
L. L. Avramov, ``Local algebra and rational homotopy,'' Homotopie 
Alg\'ebrique et Alg\`ebre locale (J.M. Lemaire and J.-C. Thomas, 
eds.), Ast\'erisque 113-114 (1984), 15-43.

\bibitem{Avr}
\rule{2cm}{.01cm}, ``Locally complete intersection homomorphisms and a
conjecture of Quillen on the vanishing of cotangent homology,''
Annals of Math. (2) 150 (1999), 455-487. 

\bibitem{AF}
L. L. Avramov and H.-B. Foxby, ``Locally Gorenstein 
homomorphisms,'' American J Math. 114 (1992), 1007-1047.  

\bibitem{AF2}
\rule{2cm}{.01cm}, ``Grothendieck's localization problem,'' Contemp. Math. 
159 (1994), 1-13.  

\bibitem{AFHal}
L. L. Avramov, H.-B. Foxby, and S. Halperin, ``Differential graded homological 
algebra,'' in preparation.

\bibitem{AFH}
L. L. Avramov, H.-B. Foxby, and B. Herzog, ``Structure of local 
homomorphisms,'' J. Algebra 164 (1994), 124-145.  

\bibitem{AFL}
L. L. Avramov, H.-B. Foxby, and J. Lescot, ``Bass series of ring
homomorphisms of finite flat dimension,'' Trans. A.M.S. 335 (2) (1993), 
497-523.  

\bibitem{AH}
L. L. Avramov and S. Halperin, ``On the non-vanishing of cotangent 
cohomology,'' Comment. Math. Helvetici  67 (1987), 169-184.  

\bibitem{BG}
A.K. Bousfield and V.K.A.M. Gugenheim, {\em On PL De Rham theory and 
rational homotopy type,} Memoirs of the A.M.S. 179 (1976).

\bibitem{Goe}
P. Goerss, ``A Hilton-Milnor theorem for categories of simplicial
algebras,'' Amer. J. Math,. 111 (1989), 927--971.

\bibitem{Mat}
H. Matsumura, {\em Commutative ring theory}, Cambridge Studies in 
Advanced Math. 8, Cambridge University Press, 1996.

\bibitem{Mil}
H. Miller, ``The Sullivan conjecture on maps from classifying spaces,''
Annals of Math. 120 (1984), 39--87.

\bibitem{Qui1}
D. Quillen, {\em Homotopical algebra}, Lecture Notes in Mathematics 43,
Springer-Verlag, 1967.

\bibitem{Qui2}
\rule{2cm}{.01cm}, ``Rational homotopy theory,'' Annals of Math. 90 
(1968), 205-295.

\bibitem{Qui3}
\rule{2cm}{.01cm}, ``On the (co)homology of commutative rings,'' Proc.
Symp. Pure Math. 17 (1970), 65--87.

\bibitem{SS}
B. Shipley and S. Schwede, ``Equivalences of monoidal model 
categories,'' Algebraic and Geometric Topology 3 (2003), 287-334.

\bibitem{Tur1}
J. M. Turner, ``On simplicial commutative algebras with vanishing 
Andr\'e-Quillen homology,'' Invent. Math. 142 (3) (2000) 547-558.

\bibitem{Tur2}
\rule{2cm}{.01cm}, ``On simplicial commutative algebras with 
Noetherian homotopy,'' J. Pure Appl. Alg. 174 (2002) pp 207-220.

\bibitem{Wei}
C. Weibel, {\em An Introduction to Homological Algebra}, Cambridge
Studies in Advanced Mathematics 38, Cambridge University Press, 1995.

\end{thebibliography}
\end{document}